\documentclass[a4paper,12pt]{article}
\usepackage[cp1251]{inputenc}
\usepackage{amsfonts, amssymb, amsmath, amsthm, amscd}
\usepackage{graphicx}
\usepackage{cite}
\ifx\undefined \pdfgentounicode \else
\input{glyphtounicode} \pdfgentounicode=1
\fi

\author{A.A. Vasil'eva}
\title{Kolmogorov widths of the intersection of two finite-dimensional balls}
\date{}
\begin{document}

\maketitle

\newenvironment{Biblio}{%
                  \renewcommand{\refname}{\footnotesize REFERENCES}%
                  \begin{thebibliography}{99}%
                  \itemsep -0.2em}%
                  {\end{thebibliography}}

\def\inff{\mathop{\smash\inf\vphantom\sup}}
\renewcommand{\le}{\leqslant}
\renewcommand{\ge}{\geqslant}
\newcommand{\sgn}{\mathrm {sgn}\,}
\newcommand{\inter}{\mathrm {int}\,}
\newcommand{\dist}{\mathrm {dist}}
\newcommand{\supp}{\mathrm {supp}\,}
\newcommand{\R}{\mathbb{R}}
\newcommand{\Z}{\mathbb{Z}}
\newcommand{\N}{\mathbb{N}}
\newcommand{\Q}{\mathbb{Q}}
\theoremstyle{plain}
\newtheorem{Trm}{Theorem}
\newtheorem{trma}{Theorem}
\newtheorem{Def}{Definition}
\newtheorem{Cor}{Corollary}
\newtheorem{Lem}{Lemma}
\newtheorem{Rem}{Remark}

\renewcommand{\proofname}{\bf Proof}
\renewcommand{\thetrma}{\Alph{trma}}

In this paper we study the problem on estimating the Kolmogorov
widths of the set $B_{p_0}^m\cap \nu B_{p_1}^m$ in $l_q^m$.

First we recall some definitions.

Let $X$ be a normed space, and let $M\subset X$, $n\in \Z_+$. The
Kolmogorov $n$-width of the set $M$ in the space $X$ is defined as
$$
d_n(M, \, X) = \inf _{L\in {\cal L}_n(X)} \sup _{x\in M} \inf
_{y\in L} \|x-y\|;
$$
here ${\cal L}_n(X)$ is the family of subspaces in $X$ of
dimension at most $n$.

Let $1\le p<\infty$. We denote by $l_p^m$ the linear space $\R^m$
with the norm $$\|(x_1, \, \dots, \, x_m)\|_{l_p^m}=\left(\sum
\limits _{j=1}^m |x_j|^p\right)^{1/p};$$ by $B_p^m$ we denote the
unit ball in $l_p^m$. For $p=\infty$ the norm is defined as
$$\|(x_1, \, \dots, \, x_m)\|_{l_\infty^m}= \max \{|x_1|, \, \dots,
\, |x_m|\}.$$

The problem on estimating the widths $d_n(B_p^m, \, l_q^m)$ was
studied by A. Pietsch, M.I. Stesin, B.S. Kashin, E.D. Gluskin,
A.Yu. Garnaev (see, e.g., \cite{gluskin1, bib_gluskin, bib_kashin,
kashin_oct, garn_glus, pietsch1, stesin}). Here we use and
formulate these results for some particular cases. First we
introduce notation for order equalities and inequalities. Let $X$,
$Y$ be sets, $f_1$, $f_2:\ X\times Y\rightarrow \mathbb{R}_+$. We
write $f_1(x, \, y)\underset{y}{\lesssim} f_2(x, \, y)$ (or
$f_2(x, \, y)\underset{y}{\gtrsim} f_1(x, \, y)$) if for each
$y\in Y$ there is $c(y)>0$ such that $f_1(x, \, y)\le c(y)f_2(x,
\, y)$ for all $x\in X$; $f_1(x, \, y)\underset{y}{\asymp} f_2(x,
\, y)$ if $f_1(x, \, y) \underset{y}{\lesssim} f_2(x, \, y)$ and
$f_2(x, \, y)\underset{y}{\lesssim} f_1(x, \, y)$.

\begin{trma}
\label{glus} {\rm \cite{bib_gluskin}} Let $1\le p\le q<\infty$,
$0\le n\le m/2$.
\begin{enumerate}
\item Let $1\le q\le 2$. Then $d_n(B_p^m, \, l_q^m) \underset{p,q}{\asymp}
1$.

\item Let $2<q<\infty$, $\lambda_{pq} =\min \left\{1, \,
\frac{1/p-1/q}{1/2-1/q}\right\}$. Then
$$
d_n(B_p^m, \, l_q^m) \underset{p,q}{\asymp} \min \{1, \,
n^{-1/2}m^{1/q}\} ^{\lambda_{pq}}.
$$
\end{enumerate}
\end{trma}

\begin{trma}
\label{p_s} {\rm \cite{pietsch1, stesin}.} Let $1\le q\le p\le
\infty$, $0\le n\le m$. Then
$$
d_n(B_p^m, \, l_q^m) = (m-n)^{1/q-1/p}.
$$
\end{trma}

The problem on estimating the Kolmogorov widths of intersections
of finite-dimensional balls naturally appears when estimating the
Kolmogorov  widths of intersections of function classes (see,
e.g., \cite{galeev1, galeev4, vas_1}). Notice that in
\cite{galeev1} the intersection of the arbitrary family of balls
was considered, but the order estimate for the widths was obtained
only for $m=2n$. In \cite{vas_1} the problem on estimating the
Kolmogorov widths of some function classes in a weighted
$L_q$-space was studied; these classes were defined by constraints
on the weighted $L_{p_1}$-norm of the highest-order derivatives
and on the weighted $L_{p_0}$-norm of the functions. The lower
estimates for the widths were reduced to estimating the widths of
$B_\infty^m$ and $B_1^m$ in $l_q^m$; the upper estimates were
reduced to estimating the widths of $B_{p_0}^m$, $B_{p_1}^m$,
$B_2^m$ and $B_q^m$ (here the particular case of Theorem 2 from
\cite{galeev1} was applied). The problem on estimating the widths
of $B_{p_0}^m\cap \nu B_{p_1}^m$ was not studied in \cite{vas_1}.

Let $1\le p_1<p_0 \le\infty$. First we notice that $B_{p_0}^m\cap
\nu B_{p_1}^m = \nu B_{p_1}^m$ for $\nu\le 1$, and $B_{p_0}^m\cap
\nu B_{p_1}^m = B_{p_0}^m$ for $\nu\ge m^{1/p_1-1/p_0}$. Hence, it
suffices to consider the case $1< \nu < m^{1/p_1-1/p_0}$. Further,
$B_{p_1}^m \subset B_{p_0}^m\cap \nu B_{p_1}^m \subset B_{p_0}^m$
for $\nu\ge 1$. Therefore, Theorem \ref{glus} implies that for
$n\le m/2$
$$
d_n(B_{p_0}^m\cap \nu B_{p_1}^m, \, l_q^m) \underset{p,q}{\asymp}
1 \quad \text{if}\; 1\le p_1\le p_0\le q\le 2,
$$
$$
d_n(B_{p_0}^m\cap \nu B_{p_1}^m, \, l_q^m) \underset{p,q}{\asymp}
\min\{1, \, n^{-1/2}m^{1/q}\} \quad \text{if}\; 1\le p_1\le p_0\le
2<q<\infty.
$$
Now, we estimate the Kolmogorov widths for other cases. Notice
that it suffices to consider $\nu = k^{1/p_1-1/p_0}$, $k=1, \,
\dots, \, m$.

\begin{Trm}
\label{main} Let $p_0>p_1$, $m\in \N$, $n\in \Z_+$, $n\le m/2$,
$1\le k\le m$, $\nu = k^{1/p_1-1/p_0}$.
\begin{enumerate}
\item Let $1\le p_1<2<p_0<q<\infty$. Then
$$
d_n(B_{p_0}^m\cap \nu B_{p_1}^m, \, l_q^m)
\underset{p_0,p_1,q}{\asymp} \left\{ \begin{array}{l} 1, \quad
\text{for} \quad n\le m^{2/q}, \\
(n^{-1/2}m^{1/q})^{\frac{1/p_0-1/q} {1/2-1/q}}, \quad
\text{for} \quad m^{\frac 2q}\le n\le k^{1-\frac 2q}m^{\frac 2q},\\
k^{\frac 12 -\frac{1}{p_0}} n^{-\frac 12}m^{\frac 1q}, \quad
\text{for} \quad k^{1-\frac 2q}m^{\frac 2q}\le n\le m/2.
\end{array}\right.
$$
\item Let $2\le p_1<p_0<q<\infty$. Then
$$
d_n(B_{p_0}^m\cap \nu B_{p_1}^m, \, l_q^m)
\underset{p_0,p_1,q}{\asymp} \left\{ \begin{array}{l} 1, \quad
\text{for} \quad n\le m^{2/q}, \\
(n^{-1/2}m^{1/q})^{\frac{1/p_0-1/q} {1/2-1/q}}, \quad
\text{for} \quad m^{\frac 2q}\le n\le k^{1-\frac 2q}m^{\frac 2q},\\
k^{\frac{1}{p_1}-\frac{1}{p_0}} \left(n^{-\frac 12}m^{\frac
1q}\right) ^{\frac{1/p_1-1/q}{1/2-1/q}}, \quad \text{for} \quad
k^{1-\frac 2q}m^{\frac 2q}\le n\le m/2.
\end{array}\right.
$$
\item Let $2\le p_1\le q\le p_0$. Then
$$
d_n(B_{p_0}^m\cap \nu B_{p_1}^m, \, l_q^m)
\underset{p_0,p_1,q}{\asymp} \left\{ \begin{array}{l}
k^{\frac{1}{q}-\frac{1}{p_0}}, \quad
\text{for} \quad n\le k^{1-\frac 2q}m^{\frac 2q},\\
k^{\frac{1}{p_1}-\frac{1}{p_0}} \left(n^{-\frac 12}m^{\frac
1q}\right) ^{\frac{1/p_1-1/q}{1/2-1/q}}, \quad \text{for} \quad
k^{1-\frac 2q}m^{\frac 2q}\le n\le m/2.
\end{array}\right.
$$
\item Let $1\le p_1<2<q\le p_0$. Then
$$
d_n(B_{p_0}^m\cap \nu B_{p_1}^m, \, l_q^m)
\underset{p_0,p_1,q}{\asymp} \left\{ \begin{array}{l}
k^{\frac{1}{q}-\frac{1}{p_0}}, \quad
\text{for} \quad n\le k^{1-\frac 2q}m^{\frac 2q},\\
k^{\frac 12 -\frac{1}{p_0}} n^{-\frac 12}m^{\frac 1q}, \quad
\text{for} \quad k^{1-\frac 2q}m^{\frac 2q}\le n\le m/2.
\end{array}\right.
$$
\item Let $q\le 2$, $1\le p_1<q<p_0$. Then
$$
d_n(B_{p_0}^m\cap \nu B_{p_1}^m, \, l_q^m)
\underset{p_0,p_1,q}{\asymp} k^{\frac 1q-\frac{1}{p_0}}.
$$
\item Let $1\le q\le p_1<p_0\le \infty$. Then
$$
d_n(B_{p_0}^m\cap \nu B_{p_1}^m, \, l_q^m)
\underset{p_0,p_1,q}{\asymp} k^{1/p_1-1/p_0}m^{1/q-1/p_1}.
$$
\end{enumerate}
\end{Trm}

Given $k\in \{1,\, \dots, \, m\}$, we denote $\hat{x}_j=1$ for
$1\le j\le k$, $\hat{x}_j=0$ for $k+1\le j\le m$,
$$
V_k={\rm conv}\, \{(\varepsilon_1 \hat{x}_{\sigma(1)}, \, \dots,
\, \varepsilon_m \hat{x}_{\sigma(m)}):\; \varepsilon_j=\pm 1, \;
1\le j\le m, \;  \sigma \in S_m\};
$$
here $S_m$ is the permutation group on $m$ elements. Notice that
$V_1 = B_1^m$, $V_m = B_\infty^m$.

We need the following E.D. Gluskin's results.

\begin{trma}
\label{gl_q_g_2} {\rm \cite{gluskin1}} There is a non-increasing
function $a(\cdot):[2, \, \infty) \rightarrow (0, \, 1]$ with the
following property: if $q\ge 2$, $1\le k\le m$, $n\le a(q)
m^{\frac 2q}k^{1 -\frac 2q}$, then
\begin{align}
\label{kq1} d_n(V_k, \, l_q^m) \underset{q}{\gtrsim} k^{1/q}.
\end{align}
\end{trma}
Also notice that arguing as in \cite{gluskin1} for $q=2$ we obtain
the inequality
\begin{align}
\label{kq2} d_n(V_k, \, l_q^m) \ge \min _{x\in
\R}\left(k-2\sqrt{\frac{kn}{m}}x+x^2\right)^{1/2} =
k^{1/2}\left(1-\frac{n}{m}\right)^{1/2}.
\end{align}

\begin{trma}
\label{gl_q_l_2} {\rm \cite{bib_glus_3}} Let $1<q\le 2$, $n\le
m/2$. Then $d_n(V_k, \, l_q^m) \underset{q}{\gtrsim} k^{1/q}$.
\end{trma}

\renewcommand{\proofname}{\bf Proof of Theorem \ref{main}}

\begin{proof}
For $p_1\le q\le p_0$ we define $\lambda\in (0, \, 1)$ by the
equation $\frac 1q =\frac{1-\lambda}{p_1} +\frac{\lambda}{p_0}$.
Then
\begin{align}
\label{q_emb} B_{p_0}^m\cap \nu B_{p_1}^m \subset \nu^{1-\lambda}
B_q^m = k^{\frac 1q-\frac{1}{p_0}}B_q^m.
\end{align}
It follows from H\"{o}lder's inequality; it is also the particular
case of Theorem 2 from \cite{galeev1}.

If $p_1\le 2\le p_0$, we define $\tilde\lambda\in (0, \, 1)$ by
the equation $\frac 12 =\frac{1-\tilde\lambda}{p_1}
+\frac{\tilde\lambda}{p_0}$; we similarly obtain
\begin{align}
\label{2_emb} B_{p_0}^m\cap \nu B_{p_1}^m \subset
\nu^{1-\tilde\lambda} B_2^m = k^{\frac 12-\frac{1}{p_0}}B_2^m.
\end{align}

Let us prove the upper estimates for the Kolmogorov widths. In
cases 1 and 2 for $n\le k^{1-\frac 2q}m^{\frac 2q}$ we use the
inclusion $B_{p_0}^m\cap \nu B_{p_1}^m \subset B_{p_0}^m$. In
cases 3 and 4 for $n\le k^{1-\frac 2q}m^{\frac 2q}$ we apply
(\ref{q_emb}); the same inclusion we use in case 5 for $n\le m/2$.

In cases 1 and 4 for $n\ge k^{1-\frac 2q}m^{\frac 2q}$ we apply
(\ref{2_emb}). In cases 2 and 3 for $n\ge k^{1-\frac 2q}m^{\frac
2q}$ we use the inclusion $B_{p_0}^m\cap \nu B_{p_1}^m \subset \nu
B_{p_1}^m$, as well as in case 6 for $n\le m/2$.

It remains to apply Theorems \ref{glus} and \ref{p_s}.

Now we obtain the lower estimate. Notice that
\begin{align}
\label{vkl} k^{-\frac{1}{p_0}} V_k \subset B_{p_0}^m\cap \nu
B_{p_1}^m,
\end{align}
\begin{align}
\label{cube} k^{1/p_1-1/p_0}m^{-1/p_1}B_\infty^m \subset
B_{p_0}^m\cap \nu B_{p_1}^m.
\end{align}

First we consider $q>2$. Let $a(\cdot)\in (0, \, 1]$ be as in
Theorem \ref{gl_q_g_2}.

If $a(q)k^{1-\frac 2q}m^{\frac 2q} \le n \le a(q)m$, there is a
number $\tilde q\in [2, \, q]$ such that
$$n=a(q)k^{1-\frac{2}{\tilde q}}m^{\frac{2}{\tilde q}}
\le a(\tilde q)k^{1-\frac{2}{\tilde q}}m^{\frac{2}{\tilde q}}.$$
Then
$$
d_n(B_{p_0}^m\cap \nu B_{p_1}^m, \, l_q^m)
\stackrel{(\ref{vkl})}{\ge} k^{-\frac{1}{p_0}} d_n(V_k, \, l_q^m)
\stackrel{(\ref{kq1})}{\underset{q}{\gtrsim}} k^{\frac{1}{\tilde
q} -\frac{1}{p_0}}\cdot m^{\frac 1q-\frac{1}{\tilde q}}
\underset{q}{\asymp} k^{\frac 12-\frac{1}{p_0}}n^{-\frac
12}m^{\frac 1q}.
$$
If $a(q)m\le n\le m/2$, we have
$$
d_n(B_{p_0}^m\cap \nu B_{p_1}^m, \, l_q^m)
\stackrel{(\ref{vkl})}{\ge} k^{-\frac{1}{p_0}} d_n(V_k, \, l_q^m)
\ge k^{-\frac{1}{p_0}} m^{\frac 1q-\frac 12} d_n(V_k, \, l_2^m)
\stackrel{(\ref{kq2})}{\underset{p_0,p_1,q}{\gtrsim}} k^{\frac
12-\frac{1}{p_0}} m^{\frac 1q-\frac 12}.
$$
Hence, we obtain the lower estimate for $a(q)k^{1-\frac
2q}m^{\frac 2q}\le n\le m/2$ in cases 1 and 4.

Now we set $\tilde k =
\left\lceil(a(q)^{1/2}n^{-1/2}m^{1/q})^{\frac{1}{1/q-1/2}}\right\rceil$.
If $a(q)k^{1-\frac 2q}m^{\frac 2q}\le n\le a(q)m$, we have $k\le
\tilde k\le m$; this yields the inequality $\nu \tilde
k^{1/p_0-1/p_1}\le 1$ and the inclusion $\nu \tilde
k^{-1/p_1}V_{\tilde k} \subset B_{p_0}^m\cap \nu B_{p_1}^m$. Since
$n\le a(q)\tilde k^{1-\frac 2q}m^{\frac 2q}$,
$$
d_n(B_{p_0}^m\cap \nu B_{p_1}^m, \, l_q^m) \ge d_n(\nu \tilde
k^{-1/p_1}V_{\tilde k}, \, l_q^m)
\stackrel{(\ref{kq1})}{\underset{q}{\gtrsim}} \nu \tilde
k^{1/q-1/p_1} \underset{p_0,p_1,q}{\asymp}
k^{\frac{1}{p_1}-\frac{1}{p_0}}(n^{-1/2}m^{1/q})^{\frac{1/p_1-1/q}{1/2-1/q}}.
$$
If $a(q)m\le n\le m/2$, we use (\ref{cube}) and get
$$
d_n(B_{p_0}^m\cap \nu B_{p_1}^m, \, l_q^m) \ge
k^{1/p_1-1/p_0}m^{-1/p_1}d_n(B_\infty^m, \,
l_q^m)\underset{p_0,p_1,q}{\gtrsim} k^{1/p_1-1/p_0}m^{1/q-1/p_1}
$$
by Theorem \ref{p_s}. Hence, we obtain the lower estimate for
$a(q)k^{1-\frac 2q}m^{\frac 2q}\le n\le m/2$ in cases 2 and 3.

Now we consider the case $n\le a(q) k^{1-\frac 2q}m^{\frac 2q}$.

Let $\tilde k = \left\lceil(a(q)^{1/2}n^{-1/2}m^{1/q})
^{\frac{1}{1/q-1/2}}\right\rceil$. Then $\tilde k\le k$; i.e.,
$\tilde k^{1/p_1-1/p_0} \le \nu$. This yields the inclusion
$\tilde k^{-1/p_0}V_{\tilde k} \subset B_{p_0}^m\cap \nu
B_{p_1}^m$. In addition, $n\le a(q)\tilde k^{1-\frac 2q}m^{\frac
2q}$. Therefore,
$$
d_n(B_{p_0}^m\cap \nu B_{p_1}^m, \, l_q^m) \ge d_n(\tilde
k^{-1/p_0}V_{\tilde k}, \, l_q^m)
\stackrel{(\ref{kq1})}{\underset{q}{\gtrsim}} \tilde k^{1/q-1/p_0}
\underset{p_0,p_1,q}{\asymp} \min(1, \, n^{-1/2}m^{1/q})
^{\frac{1/p_0-1/q}{1/2-1/q}}.
$$
Hence, we obtain the lower estimate in cases 1, 2.

Further,
$$
d_n(B_{p_0}^m\cap \nu B_{p_1}^m, \, l_q^m)
\stackrel{(\ref{vkl})}{\ge} d_n( k^{-1/p_0}V_{ k}, \, l_q^m)
\stackrel{(\ref{kq1})}{\underset{q}{\gtrsim}} k^{1/q-1/p_0}.
$$
This yields the lower estimate in cases 3, 4.

For $q\le 2$, $p_0>q>p_1\ge 1$ the inequality $d_n(B_{p_0}^m\cap
\nu B_{p_1}^m, \, l_q^m)\underset{q}{\gtrsim} k^{1/q-1/p_0}$
follows from (\ref{vkl}) and Theorem \ref{gl_q_l_2}.

For $1\le q\le p_1<p_0\le \infty$ we apply the inclusion
(\ref{cube}) together with Theorem \ref{p_s}.
\end{proof}

\end{document}